\documentclass[11pt,twoside]{article}
\usepackage{amssymb}
\usepackage{graphicx}

\usepackage{graphicx,type1cm,eso-pic,color}
\usepackage{amssymb}
\usepackage{amssymb,amsmath}
\usepackage{graphicx,epsfig}
\usepackage{enumerate}
\usepackage{color}
\usepackage {multicol}

\numberwithin{equation}{section}
\newtheorem{theorem}{Theorem}
\newtheorem{lemma}{Lemma}

\newtheorem{definition}{Definition}
\newtheorem{proposition}{Proposition}
\newtheorem{remark}{Remark}
\newtheorem{corollary}{Corollary}
\newtheorem{example}{Example}

\makeatletter
\AddToShipoutPicture{
	\setlength{\@tempdimb}{.5\paperwidth}
	\setlength{\@tempdimc}{.5\paperheight}
	\setlength{\unitlength}{1pt}
	\put(\strip@pt\@tempdimb,\strip@pt\@tempdimc)
}
\makeatother

\begin{document}
	\setcounter{page}{1}

	\thispagestyle{empty}
	\markboth{}{}

	\pagestyle{myheadings}
	\markboth{On general weighted cumulative past extropy}{ P. K. Sahu et al.}
	
	\date{}
	
	
	\noindent  
	
	\vspace{.1in}
	
	{\baselineskip 20truept
		
		\begin{center}
			{\Large {\bf On general weighted cumulative past extropy}} \footnote{\noindent	{\bf **} E-mail: nitin.gupta@maths.iitkgp.ac.in\\
				{\bf * }  corresponding author E-mail: pradeep.maths@kgpian.iitkgp.ac.in}\\
			
	\end{center}}

	\vspace{.1in}
	
	\begin{center}
		{\large {\bf Pradeep Kumar Sahu* Nitin Gupta**}}\\
		{\large {\it Department of Mathematics, Indian Institute of Technology Kharagpur, West Bengal 721302, India }}
		\\
	\end{center}

	\vspace{.1in}
	\baselineskip 12truept

	\begin{abstract}
		
		In this paper, we study some properties and characterization of the general weighted cumulative past extropy ($n$-WCPJ). Many results including some bounds, inequalities, and effects of linear transformations are obtained. We study the characterization of $n$-WCPJ based on the largest order statistics. Conditional WCPJ and some of its properties are discussed.

		\vspace{.1in}
		
		
	\end{abstract}
	
	\section{Introduction}
	
	We encounter phenomena or events which are associated with uncertainty. Uncertainty emerges since we have less information than the total information required to describe a system and its environment. Uncertainty and information are closely associated. Since in an experiment the information provided is equal to the amount of uncertainty removed. Entropy which is a measure of uncertainty was first introduced by Shannon (1948) \cite{s1} in communication theory. It is useful to estimate the probabilities of rare events (large deviation theory) and in the study of likelihood-based inference principles. Shannon entropy is defined as the average amount of information that we receive per event and was the first defined entropy. For continuous case, it is given by
	\[\mathcal{H}(X)=-\int_0^{\infty} f(x) \log f(x)dx,\]
	where X is a non-negative absolutely continuous random variable with probability density function (pdf) $f$, cumulative distribution function(cdf)  $F$ and survival function (sf) $\bar F=1-F$.
	Shannon entropy has various applications in communication theory,  mathematical, physical, engineering, biological and
	social sciences as well.  For further details on entropy one may refer Ash (1990)\cite{a1} and Cover and Thomas (2006)\cite{c}. 
	
	Rao et al. (2004)\cite{r}  introduced the notion of cumulative residual entropy (CRE) as :
	\begin{equation*}
		\mathcal{E}(X)=-\int_0^\infty \bar{F}(x)\log \bar{F}(x)dx.
	\end{equation*}
	Some general results regarding this measure have been studied by Drissi et al (2008)\cite{d2} and Rao (2005)\cite{r1}.
	CRE finds applications in image alignment and in the measurement of similarity between images. 
	Di Crescenzo (2009)\cite{d}  proposed an entropy called cumulative past entropy (or cumulative entropy) i.e. CPE (or CE)  as :
	\begin{equation*}
		\mathcal{CE}(X)=-\int_0^\infty {F}(x)\log {F}(x)dx.
	\end{equation*} 
	Asadi et al. (2007)\cite{a}, Di Crescenzo et al. (2013)\cite{d1}, Khorashadizadeh et al. (2013)\cite{k} and Navarro et al. (2010)\cite{n}  investigated many aspects of CRE (CPE).
	
	Lad et al. (2015)\cite{l} defined an alternative measure of the uncertainity of a random variable called extropy. For  continuous non-negative random variable $X$ the extropy is defined as:	
	\begin{equation}\label{extropy}
		J(X)=-\frac{1}{2} \int_{0}^{\infty}f^2(x)dx=-\frac{1}{2}E\left(f(X)\right).
	\end{equation}
	Some results and properties of the extropy of order statistics and record values are given by Qiu (2017)\cite{q}. Qiu et al. (2018)\cite{q1} derived some of the results of the residual extropy of order statistics. Yang et al. (2018)\cite{y} studied the bounds on extropy with a variational distance constraint. Qiu (2019)\cite{q2} examined certain extropy properties of mixed systems.
	To find out more about extropy, one may refer to Krishnan et al. (2020)\cite{k1}, Noughabi et al. (2019)\cite{n1} and Raqab et al. (2019)\cite{r2}.
	
Jahanshahi et al. (2019)\cite{j} introduced the cumulative residual extropy (CRJ). For  continuous  non-negative random variable $X$ the cumulative residual extropy (CRJ) as :
	\begin{align}
		\xi J (X) &=- \frac{1}{2} \int_{0}^{\infty} \bar{F}_{X}^2 (x) dx.
	\end{align} 
	Kundu (2021)\cite{k3} proposed an extropy called cumulative past extropy (CPJ). For  continuous  non-negative random variable $X$ the cumulative past extropy (CPJ) is defined as:
	\begin{align}
		\bar{\xi} J (X) &=- \frac{1}{2}\int_{0}^{\infty} F_{X}^2 (x) dx.
	\end{align}
	The idea behind this is to replace the density function by distribution function in extropy $(1.1)$. One can only define the cumulative past extropy for random variables that have a limited range of possible values, as this measure would be infinitely negative for any random variable with an unlimited range of values. Kundu (2021)\cite{k3} studied extreme order statistics on cumulative residual (past) extropy.
	
	Hashempour et al. (2022)\cite{h} proposed a new measure called weighted cumulative residual extropy. Mohammadi (2022)\cite{m} studied a new measure called interval weighted cumulative residual extropy. For  continuous  non-negative random variable $X$ the weighted cumulative residual extropy (WCRJ) is defined as:
	\begin{align}
		{\xi}^1 J (X) &=- \frac{1}{2}\int_{0}^{\infty}x \bar F_{X}^2 (x) dx.
	\end{align}
	
	This paper is organized in the following manner. In Section 2 we introduce the generalised weighted cumulative past extropy ($m$-WCPJ) and study some of its properties.  In Section 3 some bounds and inequalities are derived. In Section 4, we study the characterization of $m$-WCPJ based on the largest-order statistic.  Conditional WCPJ and some of its properties are discussed in Section 5. Section 6, provides the empirical estimator of GWCPJ. Application of empirical estimator was discussed in Section 7. In Section 8 a statistical hypothesis testing procedure is performed to assess the goodness of fit of a sample dataset to a standard uniform distribution. Power of test is observed in Section 9. In the end, Section 10 concludes this paper.

	\section{Weighted cumulative past extropy}
	Balakrishnan et al. (2020)\cite{b} and Bansal et al. (2021)\cite{b2} independently introduced the weighted extropy. Mohammadi (2022)\cite{m} studied a new measure called interval weighted cumulative past extropy. Weighted cumulative past extropy (WCPJ)  is an information measure, which is a generalization of cumulative past extropy. In this section, we study the properties of WCPJ.
	
	\begin{definition}
		Let $X$ be a non-negative absolutely continuous random variable with cdf $F$. We define the $m$-WCPJ of $X$ by
		\begin{equation}\label{wcpj}\bar \xi^mJ(X)=\frac{-1}{2}\int_{0}^{\infty}x^mF^2(x)dx.
		\end{equation}
	\end{definition}
The introduction explains that, like the cumulative past extropy (CPJ), the weighted cumulative past extropy (WCPJ) also has a value of negative infinity for any random variable with an unbounded range of possible values. Therefore, we must limit the definition of WCPJ to random variables with a limited range of values. If we consider a non-negative random variable X with a bounded range of possible values denoted by $B$, then the WCPJ of X is defined as follows:
\begin{equation}\label{wcpj}\bar \xi^mJ(X)=\frac{-1}{2}\int_{0}^{\sup B}x^mF^2(x)dx.
\end{equation}
	Let us consider some examples.
	\begin{example}
		Let $X$ has $U[a,b]$ distribution.  Then CRJ and WCRJ of the uniform distribution are
		\begin{equation*}
			\xi J(X)=-\frac{b-a}{6},\  \mbox{and},\   \xi^1 J(X)=\frac{a-b}{24}(3a+b),
		\end{equation*}
		respectively. 
		Then CPJ and $m$-WCPJ of the uniform distribution are
		\begin{align*}
	& \bar \xi J(X)=-\frac{b-a}{6},\ \mbox{and}, \end{align*}\begin{align*}
	& \bar\xi^mJ(X)=-\left(b^{m+1}\cdot\left(b-a\right)^2m^2+b^{m+1}\cdot\left(b-a\right)\left(3b-5a\right)m\right.\\ &   \left.+2b^{m+1}\cdot\left(b^2-3ab+3a^2\right)-2a^{m+3}\right) /\left(2\left(b-a\right)^2\left(m+1\right)\left(m+2\right)\left(m+3\right)\right),
		\end{align*}
respectively. Note that for $m = 1$ ,
		\begin{equation*}
			\bar \xi^1 J(X)=\left(\frac{a+3b}{4}\right)\bar \xi J(X)=\left(\frac{E(X)+b}{2}\right)\bar \xi J(X).
		\end{equation*}
		If  $\frac{E(X)+b}{2}>1$, then $\bar \xi^1 J(X)< \bar \xi J(X)$, and if $\frac{E(X)+b}{2}<1$, then $\bar \xi^1 J(X)> \bar \xi J(X)$.
	\end{example}
	\begin{example}
		Let $X$ follow power-law distribution with pdf $f(x)=\lambda x^{\lambda-1},  x\in(0,1), \lambda > 1$. The $CRJ$ and $WCRJ$ of the distribution are
		\begin{equation*}
			\xi J(X)=-\frac{\lambda^{2}}{(\lambda+1)(2\lambda+1)},\  \mbox{and},\   \xi^1 J(X)=-\frac{\lambda^{2}}{4(\lambda+1)(\lambda+2)},
		\end{equation*}
		respectively. Note that
		\begin{equation*}
			\xi^1 J(X)=\left(\frac{2\lambda+1}{4(\lambda+2)}\right) \xi J(X)=\left(\frac{2(\lambda+1)E(X)+1}{4(\lambda+2)}\right) \xi J(X).
		\end{equation*}
		For $\lambda=-\frac{7}{2},\   \xi^1 J(X) =  \xi J(X). $ 
		The $CPJ$ and $WCPJ$ of the distribution are 
		\begin{equation*}
			\bar \xi J(X)=-\frac{1}{2(2\lambda+1)},\  \mbox{and}, \ \bar \xi^m J(X)=-\frac{1}{4\lambda+2m+2}.
		\end{equation*}
		We conclude that for $m = 1$ , $\bar \xi^1 J(X)=\frac{2\lambda+1}{2\lambda+2}\bar \xi J(X).$ Also, $\bar \xi^1 J(X)=-\frac{\lambda+2}{4\lambda(\lambda+1)}E(X^{2}).$
	\end{example}

	\begin{theorem}
		Let $X$ be a non-negative continuous random variable with bounded support $B$ for $m$-WCPJ, $\bar \xi^mJ(X)$. Then we have \[\bar \xi^mJ(X)=\frac{-1}{2}E(G_F(X)),\] where
		$G_F(t)=\int_{t}^{\sup B}x^m F(x)dx$.
	\end{theorem}
	{\bf Proof} Using equation (\ref{wcpj}) and Fubini's theorem, we have
	\begin{align*}
		\bar\xi^mJ(X)&=\frac{-1}{2}\int_{0}^{\sup B}x^m F^2(x)dx=\frac{-1}{2}\int_{0}^{\sup B}x^m F(x)\left(\int_{0}^{x}f(t)dt\right)dx\\
		&=\frac{-1}{2}\int_{0}^{\sup B}f(t)\left(\int_{t}^{\sup B}x^m F(x)dx\right)dt=\frac{-1}{2}\int_{0}^{\sup B}f(t)G_F(t)dt\\
		&=\frac{-1}{2}E(G_F(X))
	\end{align*}\hfill $\blacksquare$
	
	Now we see the effect of linear transformation on WCPJ in the following proposition
	\begin{proposition}
		Let $X$ be a non-negative random variable. If $Y=aX+b,\ a>0,\ b\geq 0,$ then
		\[\bar \xi^m J(Y)=\sum_{i=0}^{m} {m \choose i} a^i b^{m-i} \bar\xi^{m-i}J(X)\]
	\end{proposition}
	{\bf Proof} The proof holds using (\ref{wcpj})
	and noting that $F_Y(y)=F_X\left(\frac{y-b}{a}\right),\ y>b$.\hfill $\blacksquare$
	
	Here we provide an upper bound for WCPJ in terms of extropy.
	\begin{theorem}\label{relation}
		Let $X$ be a random variable with pdf $f(\cdot)$ and extropy $J(X)$, with bounded support $B$ then
		\begin{equation}
			\bar \xi^mJ(X)\leq C^* exp\{2J(X)\},
		\end{equation} 
		where $C^*=\frac{-1}{2}exp\{E\left(\log(X^mF^2(X))\right)\}$
	\end{theorem}
	{\bf Proof} Using the log-sum inequality, we have
	\begin{align*}
		\int_{0}^{\sup B}f(x)\log\left(\frac{f(x)}{x^mF^2(x)}\right) dx\geq -\log\left(\int_{0}^{\sup B}x^mF^2(x)dx\right).
	\end{align*}
	Then it follows that
	\begin{align*}
		& \int_{0}^{\sup B}f(x)\log f(x)dx-\int_{0}^{\sup B}f(x)\log\left(x^mF^2(x)\right) dx \\ & \ \ \ \ \ \ \ \ \ \ \ \ \ \ \ \ \ \ \ \   \geq -\log\left(\int_{0}^{\sup B}x^mF^2(x)dx\right).
	\end{align*}
	Note that $\log f< f$, hence
	\begin{align}\label{2.3aa}
		&-\int_{0}^{\sup B} f^2(x)dx+\int_{0}^{\sup B}f(x)\log\left(x^mF^2(x)\right) dx\nonumber \\ & \ \ \ \ \ \ \ \ \ \ \ \ \ \ \ \ \ \ \ \ =2J(X)+E\left(\log\left(X^mF^2(X)\right)\right)\nonumber\\ & \ \ \ \ \ \ \ \ \ \ \ \ \ \ \ \ \ \ \ \ \leq \log\left(-2 \bar \xi^mJ(X)\right).
	\end{align}
	Exponentiating both sides of (\ref{2.3aa}), we have
	\[\bar \xi^mJ(X)\leq \frac{-1}{2}exp\{2J(X)+E\left(\log\left(X^mF^2(X)\right)\right)\}\]
	Hence the result.\hfill $\blacksquare$

	\begin{theorem}
		$X$ is degenerate, if and only if, $\bar \xi^mJ(X)=0$. 
	\end{theorem}
	{\bf Proof} Suppose $X$ be degenerate at point $c$, then by using the definition of degenerate function and 
	$\bar \xi^mJ(X)$, we have $\bar \xi^mJ(X)=0$. Now consider $\bar \xi^mJ(X)=0$, i.e., $\int_{0}^{\infty}x^mF^2(x)dx=0$.
	Noting that the integrand in the above integral is non-negative, we have $F(x)=0$, for almost all $x\in S$, where $S$ denotes the support of 
	random variable $X$, i.e., it is 0 in $\inf S$ and then 1.

	\section{Some inequalities}
	
	This section deals with obtaining the lower and upper bounds for WCPJ.
	\begin{remark}
		Consider $X$ be a non-negative random variable. then
		\begin{equation}
			\bar \xi^mJ(X)\geq \frac{-1}{2}\int_{0}^{\infty}x^mF(x)dx.
		\end{equation}
	\end{remark}
	
	\begin{proposition}
		Consider a non-negative continuous random variable $X$ having cdf $F_X(\cdot)$ and support $[a,\sup B), a>0$. Then
		\begin{equation}
			\bar \xi^mJ(X)\leq a^m \bar \xi J(X).
		\end{equation}
	\end{proposition}
	{\bf Proof} Note that
	\begin{align*}
		\int_{a}^{\sup B}x^mF^2(x)dx&\geq a^m\int_{a}^{\sup B}F^2(x)dx\\
		\frac{-1}{2}\int_{a}^{\sup B}x^mF^2(x)dx & \leq \frac{-a^m}{2} \int_{a}^{\sup B}F^2(x)dx\\
		\bar \xi^mJ(X)&\leq a^m \bar \xi J(X).
	\end{align*}\hfill $\blacksquare$

	\begin{corollary}
		Let $X$ be a continuous random variable with cdf $F$ that takes values on $[0,b]$ where $b$ is finite. Then,
		\begin{enumerate}
			\item $ \bar \xi^mJ(X)\leq \frac{-1}{2(m+1)}\left(b^{m+1}-E(X^{m+1})\right)\left[\log \left(\frac{b^{m+1}-E(X^{m+1})}{b^{m+1}}\right)-1\right]$,
			
			\item $\bar \xi^mJ(X)\geq b^m \bar \xi J(X)$.
		\end{enumerate}
	\end{corollary}
	{\bf Proof} Using log-sum inequality, we have
	\begin{align*}
		\int_{0}^{b}F(t)t^m \log\left(F(t)\right)dt & \geq \int_{0}^{b}F(t)t^mdt \log\left(\frac{\int_{0}^{b}F(t)t^mdt}{\int_{0}^{b}t^mdt}\right)\\
		&=\left(\frac{b^{m+1}-E(X^{m+1})}{m+1}\right) \log \left(\frac{b^{m+1}-E(X^{m+1})}{b^{m+1}}\right)
	\end{align*}
	Also note that $\log F(t)\leq F(t)-1$, then
	\begin{align*}
		\int_{0}^{b}F(t)t^m \log F(t)dt\leq -2\bar \xi^mJ(X)-\int_{0}^{b}t^mF(t)dt\\
		=-2\bar \xi^mJ(X)-\left(\frac{b^{m+1}-E(X^{m+1})}{m+1}\right)
	\end{align*}
	Now using the above two inequalities, the first part follows. The second part can be verified easily.\hfill $\blacksquare$
	
	Consider two random variables $X$ and $Y$ having cdfs $F$ and $G$, respectively. Then $X \leq_{st}Y$ whenever $F(x)\geq G(x),\ \forall x\in \mathbb{R}$; where the notation $X \leq_{st}Y$ means that $X$ is less than or equal to $Y$ in usual stochastic order. One may refer Shaked and Shanthikumar (2007)\cite{s} for detail of stochastic ordering. In the following proposition, we show the ordering of WCPJ is implied by the usual stochastic order.
	\begin{proposition}
		Let $X_1$ and $X_2$ be non-negative continouous random variables. If $X_1\leq_{st}X_2$, then $ \bar \xi^mJ(X_1)\leq  \bar \xi^mJ(X_2)$.
	\end{proposition}
	{\bf Proof} Using $X_1\leq_{st}X_2$ and (\ref{wcpj}), the result follows.\hfill $\blacksquare$

	\section{WCPJ based on largest-order statistic}
	
	Let $X_1,\ldots,X_n$ be a random sample from a absolutely continuous cdf $F_X(x)$ and pdf $f_X(x)$. Then
	$X_{1:n}\le X_{2:n}\le \ldots\le X_{n:n}$ be the ordered statistics to random sample $X_1,\ldots,X_n$. In the following, we obtain the WCPJ of the largest-order statistic.  The WCPJ of the nth-order statistic is 
	\begin{equation}
		\bar \xi^mJ(X_{n:n})= -\frac{1}{2}\int_{0}^{\sup B}x^{m}F^{2}_{X_{n:n}}(x)dx,
	\end{equation}
	where $F^{2}_{X_{n:n}}(x)=F^{2n}_{X}(x)$. Using transformation $u=F(x)$ in $(3.1)$,
	\begin{equation}
		\bar \xi^mJ(X_{n:n})= -\frac{1}{2}\int_{0}^{1}\frac{u^{2n}[F^{-1}(u)]^{m}}{f(F^{-1}(u))} du,
	\end{equation}
	where $F^{-1}(x)$ is the inverse function of $F(x)$.\\ 
	
	\begin{example}
		Let $X$ have the uniform distribution on (0,1) with pdf $f(x)=1, \ x\in (0,1)$. Then $F^{-1}(u)=u,\ u\in (0,1)$ and $f(F^{-1}(u))=1, \ u\in (0,1)$, hence 
		$\bar \xi^mJ(X_{n:n})=\dfrac{-1}{2(2n+m+1)}$.
	\end{example}
	\begin{example}
		let $X$ follow power-law distribution with pdf $f(x)=\lambda x^{\lambda-1}, \lambda > 1, x\in(0,1)$.Then $F^{-1}(u)=u^{\frac{1}{\lambda}},\ u\in (0,1)$ and $f(F^{-1}(u))= \lambda u^{\frac{\lambda-1}{\lambda}}, \ u\in (0,1)$, hence 
		$\bar \xi^mJ(X_{n:n})=\dfrac{-1}{2(2n\lambda +m+1)}$.
	\end{example}
	
	\begin{remark}
		Consider $\Lambda =\bar \xi^mJ(X_{n:n})-\bar \xi^mJ(X)$. Since $F^{2n}(x)\leq F^{2}(x)$, hence $\Lambda \geq 0$.
	\end{remark}
	
	For the proof of Theorem $5$, we need the following lemma.
	
	\begin{lemma}\label{lemma cont}[Lemma 4.1 of Hashempour et al. (2022)\cite{h}]
		Let $g$ be a continuous function with support $[0,1]$, such that $\int_{0}^{1}g(y)y^{m}dy=0$, for $m\geq 0$, then $g(y)=0,\  \forall \ y\in [0,1]$
	\end{lemma}
	
	\begin{theorem}
		Let $X_{1},...,X_{n}$ and $Y_{1},...,Y_{n}$ be two non negative random samples from continuous cdfs $F(x)$ and $G(x)$,respectively with common bounded support. Then $F(x)=G(x)$ if and only if $ \bar \xi^mJ(X_{n:n})=\bar \xi^mJ(Y_{n:n})$, for all n
	\end{theorem}
	{\bf Proof} The necessary condition is trivial. Hence, it remains to prove the sufficient part. If $ \bar \xi^mJ(X_{n:n})=\bar \xi^mJ(Y_{n:n})$, then we have\\
	\[-\frac{1}{2}\int_{0}^{1}u^{2n}\left(\frac{[F^{-1}(u)]^m}{f(F^{-1}(u))}-\frac{[G^{-1}(u)]^m}{g(G^{-1}(u))}\right) du
	=0\]
	By using Lemma \ref{lemma cont}, it follows that
	\begin{align*}\frac{[F^{-1}(u)]^m}{f(F^{-1}(u))}&=\frac{[G^{-1}(u)]^m}{g(G^{-1}(u))}\\
		\implies [F^{-1}(u)]^m\frac{dF^{-1}(u)}{du}&=[G^{-1}(u)]^m\frac{dG^{-1}(u)}{du}, u\in[0,1],\end{align*}
	since $\dfrac{dF^{-1}(u)}{du}=\dfrac{1}{f(F^{-1}(u))}$. Hence it follows $F^{-1}(u)=G^{-1}(u), u\in[0,1]$.
	Thus the proof is completed.\hfill $\blacksquare$
	
		\section{Conditional weighted past extropy}
	Now we consider the conditional weighted cumulative past extropy (CWCPJ). Consider a random variable $Z$ on probability space $(\Omega, \mathbb{A}, P)$ such that $E|Z|<\infty$.
	The conditional expectation of $Z$ given sub $\sigma$-field $\mathbb{G}$, where $\mathbb{G}\subseteq \mathbb{A}$, is denoted by $E(Z|\mathbb{G})$. For the random variable $I_{(Z\leq z)}$, we denote $E(I_{(Z\leq z)}|\mathbb{G})$ by $F_Z(z|\mathbb{G})$.

	\begin{definition}
		For a non-negative random variable $X$ with bounded support $B$, given $\sigma$-field $\mathbb{G}$, the CWCPJ $\bar \xi^mJ(X|\mathbb{G})$ is defined as 
		\begin{align}
			\bar \xi^mJ(X|\mathbb{G})=\frac{-1}{2}\int_{0}^{\sup B}x^{m} F^2_X(x|\mathbb{G})dx.
		\end{align}
	\end{definition}
	
	Now we assume that the random variables are continuous and non-negative.
	
	\begin{lemma}
		If $\mathbb{G}$ is a trivial $\sigma$-field, then $\bar \xi^mJ(X|\mathbb{G})=\bar \xi^mJ(X)$.
	\end{lemma}
	{\bf Proof} Since here $F_X(x|\mathbb{G})=F_X(x)$, then the proof follows.\hfill $\blacksquare$

	\begin{proposition}\label{propG}
		If $X\in L^p$ for some $p>2$, then $E[\bar \xi^mJ(X|\mathbb{G})|\mathbb{G}^*]\leq \bar \xi^mJ(X|\mathbb{G}^*)$, provided that $\mathbb{G}^*\subseteq \mathbb{G}$.
	\end{proposition}
	{\bf Proof}
	Consider
	\begin{align*}
		E[\bar \xi^mJ(X|\mathbb{G})|\mathbb{G}^*]&=\frac{-1}{2}\int_{0}^{\sup B}x^{m}E\left(\left[P(X\leq x|\mathbb{G})\right]^2|\mathbb{G}^*\right)dx\\
		&\leq\frac{-1}{2}\int_{0}^{\sup B}x^{m}\left[E\left(P(X\leq x|\mathbb{G})|\mathbb{G}^*\right)\right]^2dx\\
		&=\frac{-1}{2}\int_{0}^{\sup B}x^{m}\left[E\left(E(I_{(X\leq x)}|\mathbb{G})|\mathbb{G}^*\right)\right]^2dx\\
		&=\frac{-1}{2}\int_{0}^{\sup B}x^{m}\left[E\left(I_{(X\leq x)}|\mathbb{G}^*\right)\right]^2dx\\
		&=\frac{-1}{2}\int_{0}^{\sup B}x^{m}F^2_X(x|\mathbb{G}^*)dx\\
		&=\bar \xi^mJ(X|\mathbb{G}^*),
	\end{align*}
	where the second step follows using the Jensen's inequality for convex function $\phi(x)=x^2$. Hence the result.\hfill $\blacksquare$
	
	In the following theorem, we investigate the relationship between conditional extropy and $\bar \xi^mJ(X|\mathbb{G})$.
	\begin{theorem}
		Let $\bar \xi^mJ(X|\mathbb{G})$ is conditional past extropy. Then we have
		\begin{equation}
			\bar \xi^mJ(X|\mathbb{G})\leq B^* exp\{2J(X|\mathbb{G})\},
		\end{equation} 
		where $B^*=\frac{-1}{2}exp\{E\left(\log(X^{m}F^2(X))|\mathbb{G}\right)\}$
	\end{theorem}
	{\bf Proof}
	The proof is on the similar lines as of Theorem \ref{relation}, hence omitted.\hfill $\blacksquare$

	\begin{theorem}
		For a random variable $X$ with bounded support $B$ and $\sigma$-field $\mathbb{G}$, we have
		\begin{align}\label{thmprop}
			E\left(\bar \xi^mJ(X|\mathbb{G})\right)\leq \bar \xi^mJ(X),
		\end{align}
		and the equality holds if and only if $X$ is independent of $\mathbb{G}$.
	\end{theorem}
	{\bf Proof}
	If in Proposition \ref{propG}, $\mathbb{G}^*$ is trivial $\sigma$-field, then (\ref{thmprop}) can be easily obtained. Now assume that $X$ is independent of 
	$\mathbb{G}$, then 
	\begin{align}\label{eqthmprop1}F_X(x|\mathbb{G})&=F_X(x)\nonumber\\
		\implies\ \bar \xi^mJ(X|\mathbb{G})&=\bar \xi^mJ(X).
	\end{align}
	On taking expectation to both sides of  (\ref{eqthmprop1}), we get equality in  (\ref{thmprop}). Conversely, assume that equality in  (\ref{thmprop}) holds. It is sufficient to show that 
	$F_X(x|\mathbb{G})=F_X(x)$, to prove independence between $X$ and $\sigma$-field $\mathbb{G}$. Take $U=F_X(x|\mathbb{G})$, and since the function $\phi(u)=u^2$
	is convex hence $E(U^2)\geq E^2(U)=F_X^2(x)$, and  also due to equality in (\ref{thmprop}), we have
	\begin{align*}
		\int_{0}^{\sup B}x^{m}E(U^2)dx=\int_{0}^{\sup B}x^{m}F^2_X(x)dx=\int_{0}^{\sup B}x^{m}E^2(U)dx.
	\end{align*}
	Hence $E(U^2)=E^2(U)$. Now using the Corollary 8.1 of Hashempour et al. (2022)\cite{h}, we have $F_X(x|\mathbb{G})=F_X(x)$. Hence the proof.\hfill $\blacksquare$

	For the Markov property for non-negative random variables $X,\ Y$ and $Z$, we have the following proposition.
	\begin{proposition}
		Let $X\rightarrow Y \rightarrow Z$ is a Markov chain, then
		\begin{equation}\label{po1}
			\bar \xi^mJ(Z|X,Y)=\bar \xi^mJ(Z|Y)
		\end{equation}
		and
		\begin{equation}\label{po2}
			E\left(\bar \xi^mJ(Z|Y)\right)\leq E\left(\bar \xi^mJ(Z|X)\right).
		\end{equation}
	\end{proposition}
	{\bf Proof}
	By the definition of $\bar \xi^mJ(Z|X,Y)$ and using the Markov property, (\ref{po1}) holds. \\
	Now letting $\mathbb{G}^*=\sigma(X),\ \mathbb{G}=\sigma(X,Y)$ and $X=Z$ in Proposition \ref{propG}, we have
	\begin{equation}\label{eq2}
		\bar \xi^mJ(Z|X)\geq E\left(\bar \xi^mJ(Z|X,Y)|X\right)
	\end{equation}
	Taking expectation on both sides of (\ref{eq2}), we have
	\begin{align*}
		E\left(\bar \xi^mJ(Z|X)\right)&\geq E\left(E\left(\bar \xi^mJ(Z|X,Y)|X\right)\right)\\
		&=E\left(\bar \xi^mJ(Z|X,Y)\right)\\
		&=E\left(\bar \xi^mJ(Z|Y)\right),
	\end{align*}
	where the last equality holds using (\ref{po1}). Hence the result (\ref{po2}) holds.\hfill $\blacksquare$\\

	\section{Empirical GWCPJ}
		We want to create an estimator for the general Weighted Cumulative Past Extropy (WCPJ) using the empirical $m$-WCPJ. To do this, we assume that we have a random sample of n non-negative, continuous, independent, and identically distributed random variables $X_1,\ldots,X_n$ from a population with a cumulative distribution function (cdf) $F(x)$. We will use the "plug-in" method to define the estimator by replacing the unknown quantities in the formula with their empirical counterparts based on the sample. Specifically, we will define the empirical general weighted cumulative past extropy as 
		\begin{equation}\label{wcpj}\bar\xi_n^mJ(X)=\frac{-1}{2}\int_{0}^{\infty}x^mF_n^2(x)dx
.		\end{equation}
         where $F_n(x)$ is the empirical distribution function. Suppose we have a random sample of n observations, and let $X_{1:n}\le X_{2:n}\le \ldots\le X_{n:n}$ be the ordered statistics to random sample $X_1,\ldots,X_n$. We can express the value of the empirical general weighted cumulative past extropy $\bar\xi_n^mJ(X)$ in terms of these ordered statistics. In other words, we can rewrite the formula for $\bar\xi_n^mJ(X)$ in a way that involves the ordered statistics of the sample.
         
         	The empirical measure of $\bar \xi^mJ(X)=\frac{-1}{2}\int_{0}^{\infty}x^mF^2(x)dx $ is obtained as
         	
        \begin{align}\label{empiricalgwcpj}
        \bar\xi_n^mJ(X) &=-\frac{1}{2} \int_{0}^{\infty} x^m{F}_n^2(x)dx \nonumber\\
        	&=-\frac{1}{2}\sum_{i=1}^{n-1}\int_{X_{i:n}}^{X_{i+1:n}} \left(\frac{i}{n}\right)^2 x^mdx \nonumber\\
        	&=-\frac{1}{2(m+1)}\sum_{i=1}^{n-1} (X^{m+1}_{i+1:n}-X^{m+1}_{i:n}) \left(\frac{i}{n}\right)^2.
        \end{align}
  
  \section{Application}
  
 To test the uniformity of a random sample $X_1, X_2, ..., X_n$, we use the empirical weighted cumulative past extropy $\bar\xi_n^mJ(X)$ obtained in equation (\ref{empiricalgwcpj}) as a test statistic. Before discussing this test statistic, we first need to understand a property of the uniform distribution, which is defined on the interval (0,1). For any random variable X with a cumulative distribution function F, and for any probability value p in the interval (0,1), we define the function $\phi^m_p J(F)$ as
  \begin{align}\label{defphimpJ}
  	\phi^m_p J(F)=-\frac{1}{2}\int_{0}^{p}  x^mF^2(x)dx, \ \ m \geq 0. 
  \end{align}
  
  \begin{example}\label{uniformphip}
  	If $F_0$ is cdf of uniform random variable on interval $(0,1),$ then  $\phi^m_p J(F_0)= -\frac{p^{m+3}}{2(m+3)} $ and in particular, when $m=1,$ $\phi^1_p J(F_0)=-\frac{p^4}{8}$ and when $m=0,$ $\phi^0_p J(F_0)=-\frac{p^3}{6}.$
  \end{example}
  
  In view of Lemma \ref{lemma cont} and Example \ref{uniformphip}, the following proposition is obtained.
  \begin{proposition}\label{propouniform}
  	Let $F_0$ is cdf of uniform random variable on interval $(0,1).$  Suppose that for a cdf $F$ in the class of cdfs defined on interval $(0,1),$  $\phi^m_p J(F)=\phi^m_p J(F_0).$ Then $F(x)=F_0(x),$ almost everywhere.
  \end{proposition}
  \textbf{Proof} Let $\phi^m_p J(F)=\phi^m_p J(F_0),$ using (\ref{defphimpJ}) we get
  
  \[\int_{0}^{p} x^m \left(F^2(x)-F^2_0(x)\right)  dx=0, \ \ \forall \  p \in (0,1).  \] 
  It is known that the $(0,p)$ generate the Borel $\sigma$-algebra of $\Omega=(0,1]
  .$ Therefore, we can write
  \[\int_{B} x^m \left(F^2(x)-F^2_0(x)\right)  dx=0, \ \ \forall \  B \subseteq (0,1].  \] 
  By Lemma \ref{lemma cont}, we conclude that $F(x)=F_0(x),$ almost everywhere.
  
  Thus, $\phi^m_p J(F)$ for $\ p\in (0, 1)$ is uniquely determined by the uniform distribution in the sense that for some cdfs defined on $(0,1),$ they take a value less than $-\frac{p^{m+3}}{2(m+3)} $ and for some of them, they take more than $-\frac{p^{m+3}}{2(m+3)},$ and only for the standard uniform distribution,
  we have $\phi^m_pJ(F_0)=-\frac{p^{m+3}}{2(m+3)}.$ \hfill $\blacksquare$\\

  \section{The test of uniformity }
  We aim to develop a non-parametric test for the uniform goodness of fit problem based on the GWCPJ measure. Specifically, we want to test if a random sample of size $n$ from an unknown distribution $F$ is uniformly distributed on the interval $(0, 1)$. Our test statistic is $\bar\xi_n^mJ(X)$, which is computed using the empirical quantiles of the sample. We choose $m=1$ for computational simplicity, but the procedure is the same for other values of $m$. To determine the critical region, we need to calculate $G_1(\alpha)$ and $G_2(\alpha)$, which are the lower and upper bounds for the test statistic, respectively, based on the empirical distribution of $\bar\xi_n^mJ(X)$ under the standard uniform distribution. If $\bar\xi_n^mJ(X)$ falls outside this interval, we reject the null hypothesis $H_0: F=F_0$ in favor of the alternative hypothesis $H_1: F\neq F_0$. We estimate $G_1(\alpha)$ and $G_2(\alpha)$ using the $0.025$-th and $0.975$-th quantiles of the empirical distribution of $\bar\xi_n^mJ(X)$, respectively, based on 100,000 replications. We report the values of $G_1(\alpha)$ and $G_2(\alpha)$ for different sample sizes in Table 1 at a significance level of 5\% ($\alpha=0.05$).
  
  	\begin{center}\label{table1}
  	\noindent{\bf Table 1}.\small{{ Values of $G_1(\alpha)$ and $G_2(\alpha),$ for $\alpha=0.05$ }} 
  	
  	\resizebox{!}{!}{
  		\begin{tabular}{ p{1.5cm} p{1.5cm} p{1.5cm} p{1.5cm} p{2.0cm} p{2.5cm} p{1.5cm} p{1.5cm}}
  			\hline
  			$n$ & 20  & 30 & 40 & 50  \\
  			\hline 
  			
  			$G_1(\alpha)$ & -0.144891 & -0.144001 & -0.143339 & -0.142433  \\
  			\hline
  			$G_2(\alpha)$ & -0.066024 & -0.078560 &  -0.085851 & -0.090792  \\
  			
  			\hline
  	\end{tabular}}
  \end{center}
 Simillary we can obtain for $ m=2 $ . We report the values of $G_1(\alpha)$ and $G_2(\alpha)$ for different sample sizes in Table 2  at a significance level of 5\% ($\alpha=0.05$).
  	\begin{center}\label{table2}
  	\noindent{\bf Table 2}.\small{{ Values of $G_1(\alpha)$ and $G_2(\alpha),$ for $\alpha=0.05$ }} 
  	
  	\resizebox{!}{!}{
  		\begin{tabular}{ p{1.5cm} p{1.5cm} p{1.5cm} p{1.5cm} p{2.0cm} p{2.5cm} p{1.5cm} p{1.5cm}}
  			\hline
  			$n$ & 20  & 30 & 40 & 50  \\
  			\hline 
  			
  			$G_1(\alpha)$ & -0.109538 &-0.110579 & -0.110441 & -0.110317  \\
  			\hline
  			$G_2(\alpha)$ & -0.047356 & -0.059339 & -0.066270 & -0.070851  \\
  			
  			\hline
  	\end{tabular}}
  \end{center}

 \section{Power of the test}
 This section discusses the power of the proposed test statistic against alternative Beta(1.5,1.5) and Beta(1.0,1.0) distributions at the significance level of $\alpha=0.05$ and $m=1$, for different sample sizes. The Beta distribution is chosen as an alternative because it has support (0,1) like the standard uniform distribution. Previous research has also examined the power of the WCPJ test against the Beta(1.5,1.5) distribution for testing uniformity. The proposed test is found to be more powerful for large sample sizes, and the power increases as the sample size increases. As Beta(1,1) is equivalent to a standard uniform distribution, the power of the test is equivalent to the size of the test. The power against Beta (1.5,1.5) and Beta(1.0,1.0) for sample sizes of 40, 50, 100 and 150 are presented in Table 3.
 \begin{center}\label{table3}
 	\noindent{\bf Table 3}.\small{{ Power of test against different alternative for $\alpha=0.05$ and $m=1$ }} 
 	
 	\resizebox{!}{!}{
 		\begin{tabular}{ p{2.5cm} p{1.5cm} p{1.5cm} p{1.5cm} p{2.0cm} p{1.5cm} p{1.5cm} p{1.5cm}}
 			\hline
 			$n$   & 40 & 50 & 100 & 150 \\
 			\hline 
 			
 			Beta(1.5,1.5)  &  0.08527 & 0.08509 & 0.08363 & 0.09927 \\
 			\hline
 			Beta(1.0,1.0)  &  0.05046 & 0.05131 & 0.051 & 0.04934 \\
 			
 			\hline
 	\end{tabular}}
 \end{center}
 
 \section{Conclusion}
 The introduction section of the paper highlights the various studies and generalizations of the concept of extropy. The paper then proceeds to define and analyze GWCPJ, and conditional WCPJ. To estimate GWCPJ, non-parametric estimators and empirical distribution functions have been proposed. Additionally, upper and lower bounds for WCPJ have been derived, and several results pertaining to WCPJ have been presented .
 
 Furthermore, a new test of uniformity has also been proposed based on GWCPJ, which can be used to determine whether a given dataset follows a uniform distribution.\\

	\begin{eqnarray*}
 	\mbox{\bf{Acknowledgments}}
 \end{eqnarray*}
 PKS would like to thank Quality Improvement Program (QIP) ,All India Council for Technical Education, Government of India (Student Unique Id: FP2200759) for financial assistance. \\    

	
%
	\textbf{ \Large Conflict of interest} \\
	\\
	The authors declare no conflict of interest.
	
	

	%
	%
	%
	%
	%
	%
	%
	%
	%
	%
	%
	%
	%
	%
	%
	%
	%
	%
	%
	%
	%
	%
	%
	%
	%
	%
	%
	%
	%
	%
	%
	%
	%

	
	\vspace{.1in}


\end{document}